\documentclass[12pt]{article}
\usepackage{amssymb}

\def\i{\iota}

\def\l{\lambda}
\def\o{\omega}

\def\u{\upsilon}

\def\O{\Omega}

\chardef\tempcat=\the\catcode`\@
\catcode`\@=11
\def\cyracc{\def\u##1{\if \i##1\accent"24 i
     \else \accent"24 ##1\fi }}
\newfam\cyrfam

\def\bC{{\mathbb C}}
\def\bR{{\mathbb R}}
\def\bZ{{\mathbb Z}}

\newfont{\goth}{eufm10 scaled \magstep1}

\def\Sp#1{{\mathrm{Sp(#1)}}}

\def\SU#1{{\mathrm{SU(#1)}}}

\newfont{\mcal}{eusm10 scaled \magstep1}

\def\square{\kern1pt\vbox
             {\hrule height 0.6pt\hbox{\vrule width 0.6pt\hskip 3pt
  \vbox{\vskip 6pt}\hskip 3pt\vrule width 0.6pt}\hrule height 0.6pt}\kern1pt}

\def\ra{\rightarrow}

\def\vol{\mathrm{vol\;}}

\def\ol{\overline}

\newtheorem{Th}{Theorem}
\newtheorem{Prop}{Proposition}
\newtheorem{Cor}{Corollary}
\newtheorem{Lem}{Lemma}
\newtheorem{Def}{Definition}

\def\bt{\begin{Th}}
\def\et{\end{Th}}
\def\bp{\begin{Prop}}
\def\ep{\end{Prop}}
\def\bc{\begin{Cor}}
\def\ec{\end{Cor}}
\def\bl{\begin{Lem}}
\def\el{\end{Lem}}
\def\bd{\begin{Def}}
\def\ed{\end{Def}}

\def\pf{\noindent{\it Proof:\ }}
\def\qed{\hfill\square}
\def\n{\nabla}

\def\be{\begin{equation}}
\def\ee{\end{equation}}

\def\arr{\begin{array}{rlll}}
\def\ea{\end{array}}
\def\bea{\begin{eqnarray}}
\def\eea{\end{eqnarray}}
\def\bean{\begin{eqnarray*}}
\def\eean{\end{eqnarray*}}

\catcode`@=11
\@addtoreset{equation}{section}
\catcode`@=12

\def\K{K\"ahler} 

\def\skm{special K\"ahler manifold}
\def\km{K\"ahler manifold}
\def\lc{Levi Civita connection}
\def\kli{K\"ahlerian Lagrangian immersion}

\begin{document}
\begin{titlepage}
\rightline{December 12, 2001} 
\vskip 1.5 true cm
\begin{center}
{\LARGE  Special K\"ahler manifolds: a survey}
\vskip 1.0 true cm
{\large 
Vicente Cort\'es}
\vskip 0.8 true cm

Mathematisches Institut der Universit\"at Bonn\\
Beringstr. 1, 53115 Bonn, Germany\\
vicente@math.uni-bonn.de
\end{center}
\vskip 0.8 true cm

\begin{abstract}\noindent
This is a survey of recent contributions to the area of 
special K\"ahler geometry. It is based on lectures 
given at the 21st Winter School on Geometry and Physics held in 
Srni in January 2001. 
\end{abstract}

\end{titlepage}

\section{Remarkable features of \skm s}
A {\bf (pseudo-) \km\ } $(M,J,g)$ is a differentiable manifold 
endowed with a complex
structure $J$ and a (pseudo-) Riemannian metric $g$ such that\\ 
(i) $J$ is orthogonal with respect to the metric $g$, i.e.\ $J^*g = g$ and\\ 
(ii) $J$ is parallel with respect to the Levi Civita connection $D$, i.e.\  
$DJ = 0$.\\
In the following we will always allow 
pseudo-Riemannian, i.e.\ possibly indefinite metrics. The prefix 
``pseudo'' will be generally omitted. The following definition
is by now standard, see \cite{F}. 

\bd \label{skDef} A {\bf \skm\ } $(M,J,g,\n )$ is a \km\ 
$(M,J,g)$ together with a flat torsionfree connection $\n$ such that\\
(i) $\n \o = 0$, where $\o = g(\cdot , J \cdot )$ is the \K\ form and\\
(ii) $\n J$ is symmetric, i.e.\ $(\n_XJ)Y = (\n_YJ)X$ for all vector
fields $X$ and $Y$. 
\ed 

\noindent
More precisely, one should speak of {\it affine} \skm s since there is 
also a projective variant of \skm s. In fact, there is a 
class of (affine) \skm s $M$, which are called {\it conic} \skm s and 
which admit a certain $\bC^*$-action. The quotient of $M$ by that action 
can be considered as projectivisation of $M$ and is called a {\it projective} 
\skm , see \cite{ACD}. Originally \cite{dWVP}, 
in the supergravity literature, by a \skm\ one understood a 
projective \skm . This terminology was mantained in the first mathematical
papers on that subject \cite{C:TG,C:TAMS,AC:PLMS} and abandoned with the 
publication of \cite{F}.       

\noindent 
{\bf Example 1:} Let $(M,J,g)$ be a flat \km , i.e.\ the \lc\ $D$ is flat. 
Then $(M,J,g,\n = D)$ is a \skm\ and $\n J = 0$. Conversely, any
\skm\  $(M,J,g,\n )$ such that $\n J = 0$ satisfies $\n = D =$ \lc\ of the
flat K\"ahler metric $g$. This is the trivial example of a \skm . 

\noindent 
Before giving a general construction of \skm s, which yields plenty of
non-flat examples, I would like  to offer some motivation for that concept. 

\begin{itemize}
\item The notion of \skm\ was introduced by the physicists de Wit
and Van Proeyen \cite{dWVP} and has its origin in certain supersymmetric 
field theories. More precisely, {\it affine} special K\"ahler 
manifolds are exactly the allowed targets for the scalars of the 
vector multiplets of field theories with $N=2$ {\it rigid} 
supersymmetry on four-dimensional Minkowski spacetime. 
{\it Projective} \skm s  correspond to such theories
with {\it local} supersymmetry, which describe $N=2$ {\it supergravity}
coupled to vector multiplets. 
$N=2$ supergravity theories occur as low energy limits of type II
superstrings and play a prominent role in the study of 
moduli spaces of certain two-dimensional superconformal 
field theories \cite{CFG}. The structure of these moduli spaces is 
desribed as the product of a 
projective \skm\ and a quaternionic \km . Besides these strong 
physical motivations there is also a  number of rather mathematical
reasons to study \skm s. 
\item Interesting {\it moduli spaces} carry the structure of a \skm,
for example:
\begin{itemize}
\item The (Kuranishi) moduli space $M_X$ of gauged complex structures  
associated to a Calabi-Yau 3-fold $X$ is a \skm\ of complex signature 
$(1,n)$, $n = h^{2,1}(X)$. This fact can be found in the physical
literature, see e.g.\ \cite{S} and references therein. $M_X$ parametrises
pairs $(J, \vol )$, where $J$ is a complex structure and $\vol$ a 
$J$-holomorphic volume form on a given compact Calabi-Yau manifold
of complex dimension $3$. Let me 
recall that (from the Riemannian point of view)   
a {\bf Calabi-Yau n-fold} is a  Riemannian manifold 
with holonomy group $\SU n$.  More generally, the affine cone over any 
abstract variation of polarized Hodge structure of weight $3$ and 
with $h^{3,0} = 1$ is a (conic) \skm , see \cite{C:TAMS}. Such cones 
can be considered as formal moduli spaces, i.e.\ the underlying variation
of Hodge structure is not necessarily induced by the deformation of 
complex structure of some K\"ahler manifold. 
\item The moduli space of deformations of a 
compact complex Lagrangian submanifold $Y$ in a hyper-\km\ $X$ is a 
\skm\ with positive definite metric \cite{H:AJM}. A {\bf 
hyper-\km\ } 
is a  Riemannian manifold with holonomy group in $\Sp n$. Such a manifold
$X$ is automatically \K\ of complex dimension $2n$ and carries a holomorphic 
symplectic structure $\O$. A complex submanifold $Y \subset X$ of 
complex dimension $n$ is called {\bf Lagrangian} if $\i^*\O = 0$, where
$\iota :Y\ra X$ is the inclusion map.
\end{itemize} 
\item The cotangent bundle of any \skm\ carries the structure of a 
hyper-\km . This corresponds to the dimensional reduction  
of $N=2$ supersymmetric theories from four to three spacetime dimensions 
\cite{CFG}. This construction, which is called the
{\bf c-map} in rigid supersymmetry, 
is discussed, applied and ge\-ne\-ra\-lised in the 
mathematical literature \cite{C:TAMS,F,H:AJM,ACD}. For example, it is used 
in \cite{C:TAMS} to obtain a 
hyper-\K\ structure (of complex signature $(2,2n)$) on the bundle
${\cal J} \ra M_X$ of intermediate Jacobians over the  above moduli space 
$M_X$ associated to a Calabi-Yau 3-fold $X$. The fibre of the holomorphic 
bundle ${\cal J}$ over $(J,\vol ) \in M_X$ is the {\bf intermediate Jacobian} 
\[ \frac{H^3(X,\bC )}{H^{3,0}(X,J) + H^{2,1}(X,J) + H^3(X,\bZ )}\] 
of $(X,J)$. 
\item There is also a c-map in local supersymmetry, i.e. in supergravity, 
which to any projective \skm\ of (real) dimension $2n$ 
associates a quaternionic K\"ahler manifold of dimension 
$4n+4$ \cite{FS}. It corresponds to the dimensional reduction  
of $N=2$ supergravity coupled to vector multiplets from dimension four to 
three. For mathematical discussions of this deep construction, 
see \cite{H:manuscript,K}. 
\item The base of any algebraic completely integrable system
is a \skm , see \cite{DW,F}. An 
{\bf  algebraic completely integrable system} is a holomorphic 
submersion $\pi : X \ra M$ from a complex symplectic manifold $X$ to
a complex manifold $M$ with compact Lagrangian fibres and a smooth 
choice of polarisation on the fibres. This is essentially the 
inverse construction of the rigid c-map. There should also exist 
an inverse construction for the local c-map.
\item It was shown in \cite{ACD} that the notion of \skm\ has 
natural ge\-ne\-ra\-li\-sa\-tions in the absence of a 
metric:  ``special complex''  
 and ``special symplectic'' manifolds. The cotangent bundle of such 
manifolds carries interesting geometric structures which generalise the
hyper-K\"ahler structure on the cotangent bundle of a \skm .   
Special complex geometry (in the absence of a metric) 
may provide insight in physical theories for which no Lagrangian formulation
(and for that reason no target metric) is available. 
\item There is a close relation between \skm s and affine 
differential geometry discovered in \cite{BC1}. In fact, any 
simply connected \skm\ has a canonical realisation as a parabolic affine
hypersphere. This will be explained in detail in section \ref{AffSec}. 
\item Any projective \skm\ has a canonical (pseudo-) Sasakian circle 
bundle which is realised as a proper affine hypersphere \cite{BC3}. 
\item It was discovered in \cite{BC2} that special K\"ahler manifolds with 
a flat indefinite metric have a nontrivial moduli space, which is closely 
related to the moduli space of Abelian simply transitive affine 
groups of symplectic type.   
\item Homogeneous projective \skm s were classified, under various 
assumptions in \cite{dWVP:CMP,C:TG,AC:PLMS}. Under the c-map they 
give rise to homogeneous quaternionic 
K\"ahler manifolds.   If one restricts attention to the homogeneous projective 
\skm s of {\it semi-simple} 
group, then one finds a list of Hermitian symmetric spaces of 
non-compact type which shows a 
remarkable coincidence with the list of irreducible  
special holonomy groups of torsionfree symplectic connections \cite{MS}, 
as was noticed in \cite{AC:PLMS}. Finally, the classification \cite{AC:PLMS} 
may lead to the generalisation 
of recent ideas of Hitchin about special features of 
geometry in six dimensions to other dimensions \cite{H:AT}. 
\end{itemize}
\section{The construction of special K\"ahler manifolds}
In this section we will see that the equations defining special
K\"ahler manifolds are completely integrable, in the sense that
the general local solution can be obtained from a free holomorphic
potential. The discussion follows \cite{ACD} and is based on the 
extrinsic approach to special K\"ahler manifolds developped 
in \cite{C:TAMS}. For a similar discussion from the bi-Lagrangian
point of view see \cite{H:AJM}. 

The ambient data in the extrinsic approach 
are the following: 
The complex symplectic vector space $V = T^*{\bC}^n = {\bC}^{2n}$ with
canonical coordinates $(z^1, \ldots , z^n, w_1, \ldots ,w_n)$. In these
coordinates the symplectic form is 
\[ \O = \sum_{i=1}^n dz^i\wedge dw_i\, .\]
We denote by $\tau :V \ra V$ the complex conjugation with respect to
$V^{\tau} = T^*{\bR}^n = {\bR}^{2n}$. The algebraic data $(V, \O , \tau )$
induce on $V$ the Hermitian form
\[ \gamma := \sqrt{-1}\O (\cdot ,\tau \cdot ) \]
of complex signature $(n,n)$. 

Let $M$ be a connected complex manifold of complex dimension $n$. We denote
its complex structure by $J$. 

\bd A holomorphic immersion $\phi : M\ra V$ is called
{\bf K\"ahlerian} if 
$\phi^*\gamma$ is nondegenerate and it is called {\bf Lagrangian}
if $\phi^*\O = 0$. 
\ed 
A  K\"ahlerian immersion $\phi : M\ra V$ induces on $M$ the 
pseudo-Riemannian metric $g = {\rm Re}\, \phi^*\gamma$ such that
$(M,J,g)$ is a K\"ahler manifold. 

\bl Let $\phi : M\ra V$ be a \kli . Then the \K\ form 
$\o = g(\cdot , J\cdot )$ of the K\"ahler manifold $(M,J,g)$ is given by
\[ \o = 2\sum_{i=1}^nd\tilde{x}^i\wedge d\tilde{y}_i\, ,\]
where $\tilde{x}^i := {\rm Re}\, \phi^*z^i$ and $\tilde{y}_i :=
{\rm Re}\, \phi^*w_i$. 
\el 

\pf The metric $g_V := {\rm Re}\, \gamma$ is a flat K\"ahler metric
of (real) signature $(2n,2n)$ on the complex vector space $(V,J)$. Its
K\"ahler form is 
\[ \o_V := \sum (dx^i\wedge dy_i + du^i\wedge dv_i) \, ,\]
where $x^i := {\rm Re}\, z^i$, $y_i :=
{\rm Re}\, w_i$, $u^i := {\rm Im}\, z^i$ and $v_i := {\rm Im}\, w_i$. 
On the other hand, the two-form 
\[ {\rm Re}\, \O = \sum (dx^i\wedge dy_i - du^i\wedge dv_i)\]
vanishes on $M$. This shows that 
\[ \o = \phi^*\o_V = 2 \sum \phi^*(dx^i\wedge dy_i) = 2 \sum 
d\tilde{x}^i\wedge d\tilde{y}_i\, .\qed \] 

\noindent 
The lemma implies that the functions $\tilde{x}^1, \dots , \tilde{x}^n, 
\tilde{y}_1, \ldots , \tilde{y}_n$ define local coordinates near each point
of $M$. Therefore we can define a flat torsionfree connection $\n$ on $M$
by the condition $\n d\tilde{x}^i  = \n d\tilde{y}_i = 0$, 
$i = 1,  \dots , n$. Now we can formulate the following fundamental
theorem.  

\bt  Let $\phi : M\ra V$ be a \kli\ with induced geo\-me\-tric data
$(g,\n )$. Then $(M,J,g,\n )$ is a \skm . Conversely, any simply connected
\skm\  $(M,J,g,\n )$ admits a \kli\ $\phi : M\ra V$ inducing the 
data $(g,\n )$ on $M$. The \kli\ $\phi$ is unique up to an affine
transformation of $V = {\bC}^{2n}$ with linear part in $\Sp{{\bR}^{2n}}$. 
\et 

\noindent 
For the proof of that result and its projective version see \cite{ACD}, where 
analogous extrinsic characterisations are obtained also for special 
complex and special symplectic manifolds.  

The above theorem may be considered as an extrinsic reformulation of the 
intrinsic De\-fi\-ni\-tion \ref{skDef}. The important advantage of the 
extrinsic
characterisation in terms of \kli s lies in the well known 
fact that Lagrangian immersions are locally defined by a 
generating function. More precisely, any holomorphic 
Lagrangian immersion into $(V,\O )$ is 
locally defined by a holomorphic function  
$F(z^1, \ldots ,z^n)$, at least after suitable choice of canonical 
coordinates $(z^1, \ldots , z^n, w_1, \ldots ,w_n)$. In fact, such a 
function defines a Lagrangian local section $\phi = dF$ of $T^*\bC^n = V$.
It is a \kli\ if it satisfies the nondegeneracy condition 
$\det {\rm Im}\, \partial^2 F \neq 0$. Similarly, 
projective special K\"ahler manifolds are locally defined by 
a holomorphic function $F$ satisfying a nondegeneracy condition
and which in addition is homogeneous of degree $2$.

\section{Special K\"ahler manifolds as affine hyperspheres}
\label{AffSec}
The main object of affine differential geometry are hypersurfaces
in affine space $\bR^{m+1}$ with its standard connection denoted by 
$\widetilde{\n}$ and parallel volume form ${\rm vol}$. A hypersurface 
is given by an immersion 
$\varphi : M \rightarrow \bR^{m+1}$ of an $m$-dimensional connected 
manifold. We assume that $M$ admits a transversal vector field $\xi$ and 
that $m>1$. This induces
on $M$ the volume form $\nu = {\rm vol} 
(\xi , \ldots )$, a torsionfree connection
$\n$, a quadratic covariant tensor field $g$, an endomorphism field $S$ 
(shape tensor) and a 
one-form $\theta$ such that
\bean \widetilde{\n}_XY &=& \n_XY + g(X,Y)\xi\, ,\\ 
 \widetilde{\n}_X\xi &=& SX + \theta (X)\xi \, .
\eean
We will assume that $g$ is nondegenerate and, hence, is a 
pseudo-Riemannian metric on $M$. This condition does not depend
on the choice of $\xi$.  
According to Blaschke \cite{Blaschke}, 
once the orientation of $M$ is fixed, there is a unique
choice of transversal vector field $\xi$ such that
$\nu$ coincides with the metric volume form ${\rm vol}^g$ and $\n \nu = 0$. 
This particular choice of  transversal vector field is called the 
{\bf affine normal} and the corresponding geometric data $(g,\n )$ are called
the {\bf Blaschke data}. 
Notice that, for the affine normal, $\theta = 0$ and 
$S$ is computable from $(g,\n )$ (Gau{\ss} equations).  
Henceforth we use always the affine normal as transversal
vector field.  
\bd \label{sphereDef} 
The hypersurface $\varphi : M \rightarrow \bR^{m+1}$ is called
a {\bf parabolic} (or improper) {\bf hypersphere} if the affine
normal is parallel, $\widetilde{\n}\xi = 0$. It is called a 
{\bf proper hypersphere} if the lines generated by the affine normals
intersect in a point $p \in \bR^{m+1}$, which is called {\bf the centre}. 
For parabolic hyperspheres the centre is at $\infty$. 
\ed    

\noindent 
Notice that $\widetilde{\n}\xi = 0 \Leftrightarrow S = 0 
\Leftrightarrow \n$ is flat. For proper hyperspheres $S = \l {\rm Id}$, 
$\l \in \bR -\{ 0\}$.

\noindent 
The main result of \cite{BC1} is the following:
\bt \label{BCThm} Let $(M,J,g,\n )$ be a simply connected 
special K\"ahler manifold. Then 
there exists a parabolic hypersphere $\varphi : M \rightarrow \bR^{m+1}$, 
$m = \dim_{\bR} M = 2n$, with Blaschke data $(g,\n )$. 
The immersion $\varphi$ is unique up to a unimodular affine transformation
of $\bR^{m+1}$. 
\et  

\noindent 
The proof of Theorem \ref{BCThm} makes use of 
the Fundamental Theorem of affine differential geometry \cite{DNV}, which  
is the generalisation of Radon's theorem \cite{R} to higher dimensions:  

\bt \label{DNVThm} Let $(M,g,\n )$ be a simply connected oriented 
pseudo-Riemannian manifold with a torsionfree connection $\n$ such that the 
Riemannian volume form ${\rm vol}^g$ is $\n$-parallel. Then there exists
an immersion $\varphi : M\ra  \bR^{m+1}$ with Blaschke data $(g,\n )$ if 
and only if the $g$-conjugate connection $\ol{\n}$ is torsionfree and 
projectively flat. The immersion is unique up to unimodular 
affine transformations of $\bR^{m+1}$.  
\et 
Recall that the $g$-conjugate connection $\ol{\n}$ on $M$ is defined by the 
equation:
\[ Xg(Y,Z) = g(\nabla_XY,Z) + g(Y,\overline{\nabla}_XZ) \quad 
\mbox{for all vector fields} \quad X,Y,Z \, .\]  

\pf  (of Theorem \ref{BCThm})  Let $(M,J,g,\n )$ be a simply connected 
special K\"ahler manifold. For any \km\ we have ${\rm vol}^g = 
\frac{\o^n}{n!}$, $n = \dim_{\bC} M$. This implies the first
integrability condition $\n {\rm vol}^g = 0$, since $\n \o = 0$. 
The conjugate connection is torsionfree and flat. This will follow 
{}from: 
\bl \cite{BC1} Let $(M,J,g,\n )$ be a special K\"ahler manifold. 
Then $\ol{\n} = \n^J$, where the connection $\n^J$ is defined by 
\[ \n_X^JY := J\n_X(J^{-1}Y) = - J\n_X(JY) = \n_XY -J(\n_XJ)Y \]
for all vector fields $X$ and $Y$. 
\el  
{}From the formula defining $\n^J$ we see that a vector field $X$ is 
$\n$-parallel if and only if $JX$ is $\n^J$-parallel. Therefore, if 
$X_i$, $i=1, \ldots, 2n = \dim_{\bR} M$, 
is a parallel local frame for the flat connection $\n$ 
then $JX_i$, $i=1, \ldots, 2n$, 
is a $\n^J$-parallel local frame. This shows that $\ol{\n} = \n^J$ is flat. 
Similarly the torsionfreedom of $\ol{\n} = \n^J = \n -J\n J$ 
follows from that of $\n$ and the symmetry of the tensor $\n J$. 
So the assumptions of Theorem \ref{DNVThm} are satisfied and we conclude
the existence of a hypersurface $\varphi : M \rightarrow \bR^{m+1}$, 
$m = 2n$. inducing on $M$ the Blaschke data $(g,\n )$. Now the flatness
of $\n$ implies that $\varphi$ is a parabolic hypersphere. \qed 

\noindent
As applications we obtain: 
\bc \label{constrCor} 
Any holomorphic function $F$ on a simply connected open set $U \subset
\bC^n$ with $\det {\rm Im}\, \partial^2 F \neq 0$ defines a parabolic 
hypersphere of dimension $m=2n$. 
\ec  
This is a generalisation of a classical theorem of Blaschke about 
$2$-dimensional
parabolic spheres. An explicit representation formula for the parabolic 
hypersphere in terms of the holomorphic function $F$ was given in 
\cite{C3}.

\bc  \label{rigidityCor} 
Let $(M,\nabla , g)$ be a special K\"ahler manifold with (positive) 
definite metric $g$. If $g$ is complete then $\n$ is the \lc\ and $g$ is 
flat.  
\ec 

\pf This follows by combining Theorem \ref{BCThm} with the 
following classical theorem of Calabi and Pogorelov \cite{Ca,P}. \qed  

\bt
If the Blaschke metric $g$ of a parabolic affine 
hypersphere $(M,g,\n)$ is definite 
and complete, then 
$M$ is affinely congruent to the paraboloid $x^{m+1} 
= \sum_{i=1}^m (x^i)^2$ in $\bR^{m+1}$. 
In particular, $\n$ is the \lc\ and $g$ is flat. 
\et

\noindent 
Lu \cite{Lu} proved that any special K\"ahler manifold $(M,g,J,\n )$ 
with a definite and complete metric $g$ is flat without making use
of Calabi and Pogorelov's Theorem. Special \km s $(M,J,g,\n)$ with a flat 
indefinite and geodesically complete metric $g$ for which $\n$ is complete
and is {\it not} the \lc\ were constructed in \cite{BC2}.   

\noindent 
For projective \skm s $\ol{M}$ it has been established in \cite{BC3} that 
a natural circle bundle $S \ra \ol{M}$ can be canonically realised as a 
proper hypersphere. Moreover, the metric cone over $S$ is a 
conic special K\"ahler manifold $M$, which is in turn realised as a 
parabolic hypersphere in a compatible way. There are also projective  
analogues of Corollaries \ref{constrCor} and \ref{rigidityCor}.

\end{document}